\documentclass[11pt]{article}

\usepackage{amsmath,amsfonts,amsthm,amssymb,tikz,tikz-cd,url,hyperref,cleveref,old-arrows}

\pdfoutput=1

\theoremstyle{plain}
\newtheorem{thm}{Theorem}[section]

\newtheorem{cor}[thm]{Corollary}

\theoremstyle{definition}

\newtheorem{rem}[thm]{Remark}

\theoremstyle{remark}

\newcommand{\bbA}{\mathbb{A}}
\newcommand{\bbB}{\mathbb{B}}
\newcommand{\bbC}{\mathbb{C}}

\newcommand{\bbF}{\mathbb{F}}

\newcommand{\bbN}{\mathbb{N}}

\newcommand{\bbQ}{\mathbb{Q}}
\newcommand{\bbR}{\mathbb{R}}

\newcommand{\bbZ}{\mathbb{Z}}

\newcommand{\calP}{\mathcal{P}}

\newcommand{\frakp}{\mathfrak{p}}

\newcommand{\al}{\alpha}

\newcommand{\Gam}{\Gamma}

\newcommand{\Del}{\Delta}

\newcommand{\thet}{\theta}

\newcommand{\Lam}{\Lambda}
\newcommand{\sig}{\sigma}

\DeclareMathOperator{\U}{U}

\DeclareMathOperator{\SL}{SL}

\DeclareMathOperator{\SO}{SO}

\DeclareMathOperator{\PU}{PU}

\DeclareMathOperator{\Hom}{Hom}

\DeclareMathOperator{\res}{res}

\DeclareMathOperator{\im}{Im}

\newcommand{\bs}{\backslash}

\newcommand{\lra}{\longrightarrow}

\newcommand{\wh}{\widehat}
\newcommand{\wt}{\widetilde}

\newenvironment{pf}{\begin{proof}}{\end{proof}}

\newenvironment{enum}{\begin{enumerate}}{\end{enumerate}}

\allowdisplaybreaks

\usetikzlibrary{arrows}

\title{Cohomological nonvanishing for algebraic fundamental groups of ball quotients}
\author{Matthew Stover \\ \small{Temple University}\\ \small{\textsf{mstover@temple.edu}}}
\date{\today}

\begin{document}

\maketitle

%%%%%%%%%%%%%%%%%%%%
\begin{abstract}
Suppose $\Gam < \PU(n,1)$ is a cocompact arithmetic lattice of simplest type with profinite completion $\wh{\Gam}$. This paper proves there is an open subgroup $\wh{\Gam}_0 \le \wh{\Gam}$ such that $H^j(\wh{\Del}, \bbF_p)$ is nontrivial for every open subgroup ${\wh{\Del} \le \wh{\Gam}_0}$, $j \le 2n$, and sufficiently large prime $p$. If $n \ge 2$, nonvanishing is new for all $j \ge 2$. Consequently, the virtual cohomological dimension of $\wh{\Gam}$ is at least $2n$, improving the previous lower bound of $1$. The proof shows there is a profinite fundamental class for the associated ball quotient and that its canonical class is profinite modulo torsion. For congruence $\Gam$ and $j < \frac{n+1}{2}$, restriction ${H^j(\wh{\Gam}, \bbF_p) \to H^j(\Gam, \bbF_p)}$ is shown to be almost surjective in a precise sense; this is related to whether lattices in $\PU(n,1)$ are good in the sense of Serre, which is only known to hold for $n=1$.
\end{abstract}
%%%%%%%%%%%%%%%%%%%%

%%%%%%%%%%%%%%%%%%%%
\section{Introduction}\label{sec:Intro}
%%%%%%%%%%%%%%%%%%%%

A group $\Gam$ with profinite completion $\wh{\Gam}$ is called \emph{good} if the restriction map
\begin{equation}\label{eq:Good}
H^j(\wh{\Gam}, M) \lra H^j(\Gam, M)
\end{equation}
is an isomorphism for all finite $\wh{\Gam}$ modules $M$ and all degrees $j$ \cite[\S I.2.6]{Serre}. Canonical examples include finitely generated abelian groups, fundamental groups of all closed $d$-manifolds for $d \le 3$, and groups commensurable with them; see \cite[\S 7]{Reid1}. Goodness is moreover closed under direct products \cite[Prop.\ 3.4]{GJZZ} and can be well-behaved under extensions \cite[Lem.\ 7.2]{Reid1}, leading to an expanded yet still fairly limited roster of examples.

Goodness is particularly interesting when $\Gam$ is the fundamental group of a smooth complex projective variety $X$, since $\wh{\Gam}$ is then the \emph{algebraic fundamental group} of $X$. For example, if $X$ is a $K(\Gam, 1)$ then the associated scheme is an \emph{algebraic} $K(\wh{\Gam}, 1)$ (see \cite[Apx.\ A]{Stix}) if and only if $\Gam$ is good \cite[\S 7]{Lochak}. However, the only aspherical smooth complex projective varieties known to determine an algebraic $K(\wh{\Gam}, 1)$ are essentially, via the above remarks on goodness, quotients of direct products of abelian varieties and hyperbolic smooth projective curves. Expanding the catalogue of smooth projective schemes that are an algebraic $K(\wh{\Gam}, 1)$ is a difficult and interesting problem requiring a better understanding of the maps in \Cref{eq:Good} for $j \ge 2$.

The natural aspherical generalizations of hyperbolic smooth projective curves are the \emph{ball quotients} $\Gam \bs \bbB^n$, where $\Gam$ is a torsion-free cocompact lattice in the Lie group $\PU(n,1)$ and $\bbB^n$ is the unit ball in $\bbC^n$ with its Bergman metric of constant holomorphic sectional curvature $-1$. Very little is known about algebraic fundamental groups of ball quotients (see \cite{StoverProfinite, StoverFake} for some results), not to mention whether or not there is a lattice in $\PU(n,1)$ for $n \ge 2$ that is good. Toward that end, the primary aim of this paper is to prove nonvanishing results for $H^j(\wh{\Gam}, \bbF_p)$, always meaning trivial coefficients in the finite field $\bbF_p$ with cardinality a prime $p$, among the most widely-studied class of lattices, namely those that are \emph{arithmetic of simplest type}; see \Cref{sec:Theta} for a precise definition. The first main result of this paper gives nonvanishing up to the cohomological dimension of $\Gam \bs \bbB^n$.

%%%%%%%%%%%%%%%%%%%%
\begin{thm}\label{thm:Main}
Let $\Gam < \PU(n,1)$ be a cocompact arithmetic lattice of simplest type and $\wh{\Gam}$ be its profinite completion. Then there is an open subgroup $\wh{\Gam}_0 \le \wh{\Gam}$ such that the cohomology groups $H^j(\wh{\Del}, \bbF_p)$ with trivial $\bbF_p$ coefficients are nontrivial for every open subgroup ${\wh{\Del} \le \wh{\Gam}_0}$, $j \le 2n$, and sufficiently large prime $p$.
\end{thm}
%%%%%%%%%%%%%%%%%%%%

Open subgroups of $\wh{\Gam}$ are in one-to-one correspondence with profinite completions of finite index subgroups of $\Gam$ \cite[\S 4.2]{Reid1}, hence \Cref{thm:Main} can be restated in terms of passage to a finite index subgroup $\Del \le \Gam$. The case $n = 1$, as with all the results in this paper, is well-known thanks to goodness of lattices in $\PU(1,1)$. \Cref{thm:Main} is also not new in the lowest degrees. Indeed, the map $H^j(\wh{\Gam}, \bbF_p) \to H^j(\Gam, \bbF_p)$ is an isomorphism for $j = 0,1$ \cite[\S I.2.6]{Serre} and any $\Gam$ as in the statement of \Cref{thm:Main} has a finite index subgroup $\Gam_0$ with $H^1(\Gam_0, \bbZ)$ infinite \cite[Thm.\ 1]{Kazhdan}, hence all open $\wh{\Del} \le \wh{\Gam}_0$ have $H^1(\wh{\Del}, \bbF_p)$ nontrivial for all primes $p$. \Cref{thm:Main} is the first indication that $H^j(\wh{\Gam}, \bbF_p)$ can be nontrivial all the way up to degree $2n$, which is the virtual cohomological dimension of $\Gam$.

To be precise, the $p$-cohomological dimension $\mathrm{cd}_p(G)$ of a profinite group $G$ is the least $j$ so that $H^r(G, M)$ has trivial $p$-primary component for all finite $G$-modules $M$ and all $r > j$, and the cohomological dimension of $G$ is then the supremum $\mathrm{cd}(G)$ of $\mathrm{cd}_p(G)$ over all primes $p$. The virtual cohomological dimension $\mathrm{vcd}(G)$ is the infimum of the cohomological dimensions of all the open subgroups of $G$. \Cref{thm:MainTech} then leads to the following significant improvement on the previously-known bound $\mathrm{vcd}(\wh{\Gam}) \ge 1$ easily derived from the above discussion of $H^1(\Gam, \bbZ)$.

%%%%%%%%%%%%%%%%%%%%
\begin{cor}\label{cor:Main}
Let $\Gam < \PU(n,1)$ be a cocompact arithmetic lattice of simplest type and $\wh{\Gam}$ be its profinite completion. Then $\mathrm{vcd}(\wh{\Gam}) \ge 2n$.
\end{cor}
%%%%%%%%%%%%%%%%%%%%

%%%%%%%%%%%%%%%%%%%%
\begin{pf}
Suppose $\wh{\Lam} \le \wh{\Gam}$ is an open subgroup. If $\wh{\Gam}_0$ is the open subgroup of $\wh{\Gam}$ from \Cref{thm:Main}, then $H^{2n}(\wh{\Lam} \cap \wh{\Gam}_0, \bbF_p)$ is nontrivial for all but finitely many primes $p$. The $\wh{\Lam}$-module $M_p$ induced from $\bbF_p$ then has $H^{2n}(\wh{\Lam}, M_p)$ nontrivial by Shapiro's lemma \cite[\S 1.2.5]{Serre}. Thus for every open subgroup $\wh{\Lam} \le \wh{\Gam}$ there is a finite $\wh{\Lam}$-module $M_p$ with $H^{2n}(\wh{\Lam}, M_p)$ nontrivial, hence $\mathrm{vcd}(\wh{\Gam}) \ge 2n$.
\end{pf}
%%%%%%%%%%%%%%%%%%%%

Using prior joint work with Toledo \cite[Cor.\ 3.6]{StoverToledo2}, the method of proof of \Cref{thm:Main} shows that the canonical class and its exterior powers, including the fundamental class, are `profinite modulo torsion'. More precisely:

%%%%%%%%%%%%%%%%%%%%
\begin{thm}\label{thm:MainTech}
Let $\Gam < \PU(n,1)$ be a cocompact arithmetic lattice of simplest type. Then there is a torsion-free $\Gam_0 \le \Gam$ of finite index so that for every finite cover $X = \Del \bs \bbB^n$ of $\Gam_0 \bs \bbB^n$ there is a fixed torsion class
\[
\sig \in H^2(X, \bbZ) \cong H^2(\Del, \bbZ)
\]
such that the reduction modulo $p$ of the $i^{th}$ exterior power of $K_X + \sig$ satisfies
\[
0 \neq (K_X + \sig)^i \!\!\!\!\!\pmod{p} \in \im\!\left(\mathrm{res} : H^{2i}(\wh{\Del}, \bbF_p) \lra H^{2i}(\Del, \bbF_p)\right)
\]
for all $i \le n$ and sufficiently large primes $p$, where $K_X \in H^2(X, \bbZ)$ denotes the first Chern class of the canonical line bundle on $X$.
\end{thm}
%%%%%%%%%%%%%%%%%%%%

The next natural question is whether or not $H^j(\wh{\Gam}, \bbF_p)$ can be made arbitrarily large on open subgroups. The final main result of this paper is that this can indeed be achieved for sufficiently small $j$. Moreover, toward the question of whether or not $\Gam$ is good, the map to $H^j(\Gam, \bbF_p)$ is virtually surjective in a precise sense.

%%%%%%%%%%%%%%%%%%%%
\begin{thm}\label{thm:Main2}
Suppose $\Gam < \PU(n,1)$ is a cocompact congruence arithmetic lattice of simplest type with profinite completion $\wh{\Gam}$, $j < \frac{n+1}{2}$, and $\phi_1, \dots, \phi_r \in H^j(\Gam, \bbZ)$ generate a free submodule of rank $r$. Then, for all sufficiently large primes $p$, there is a finite index subgroup $\Del_r \le \Gam$ independent of $p$ so that the restrictions of each $\phi_i$ to $\Del_r$ reduce modulo $p$ to a rank $r$ subspace contained in the image of the restriction ${H^j(\wh{\Del}_r, \bbF_p) \to H^j(\Del_r, \bbF_p)}$.
\end{thm}
%%%%%%%%%%%%%%%%%%%%

%%%%%%%%%%%%%%%%%%%%
\subsubsection*{Comments on the main results}
%%%%%%%%%%%%%%%%%%%%

%%%%%%%%%%%%%%%%%%%%
\begin{rem}
It is possible that \Cref{thm:Main} and \Cref{thm:Main2} hold for all primes, not just all but finitely many. It is also reasonable to think that the torsion class $\sig$ in \Cref{thm:MainTech} can be trivial, hence the canonical class and its exterior powers are genuinely profinite. Note however that the canonical class can be a nontrivial multiple of another class in $H^2(X, \bbC)$ (see \cite{StoverToledo1, StoverToledo2} for examples where it is divisible by arbitrarily large integers), so the finite exceptional set of primes in \Cref{thm:MainTech} is necessary.
\end{rem}
%%%%%%%%%%%%%%%%%%%%

%%%%%%%%%%%%%%%%%%%%
\begin{rem}
It is unclear if analogous results hold for every lattice in $\PU(n,1)$. Underpinning this paper is the fact that much of $H^*(X, \bbC)$ is in the span to the Poincar\'e duals to the complex totally geodesic subspaces or  is virtually generated by the cup product from $H^1(X, \bbC)$. Nonvanishing of $H^1(\Gam, \bbC)$ for congruence lattices is restricted to lattices of simplest type \cite{Clozel}, and there are $\Gam \bs \bbB^n$ with no proper complex totally geodesic subvarieties expect points. For nonarithmetic lattices, there are only finitely many maximal totally geodesic subspaces \cite{BFMS2, BU} and little is known about the cup product structure on cohomology. Thus completely different tools may be required to analyze nonvanishing of $H^j(\wh{\Gam}, \bbF_p)$ for other $\Gam$.
\end{rem}
%%%%%%%%%%%%%%%%%%%%

%%%%%%%%%%%%%%%%%%%%
\begin{rem}
Lattices in higher-rank simple Lie groups with the congruence subgroup property fail to be good, for the same reason given for $\SL_3(\bbZ)$ in \cite[\S 7.3]{Reid1}. Fundamental groups of many compact (real) hyperbolic manifolds are good, for example the lattices in $\SO(n,1)$ of simplest type \cite[Thm.\ 6.5]{KRS}. Lattices in $\PU(n,1)$ in many ways lie in between these two cases, so it is unclear whether or not one expects them to be good. It may be possible using some of the very explicit constructions in the literature to analyze an example; e.g., see \cite{Hirzebruch, BHH, Livne, Kapovich, ToledoMap, DerauxForget, DiCerboStover1, DiCerboStover2, DiCerboStover3, StoverExs} for a wealth of concrete ball quotient constructions where one could perhaps analyze cohomological properties of profinite completions closely enough to determine goodness.
\end{rem}
%%%%%%%%%%%%%%%%%%%%

%%%%%%%%%%%%%%%%%%%%
\begin{rem}
For $n\ge 2$ and $\Gam$ a lattice in $\PU(n,1)$, no upper bound is known for $\mathrm{vcd}(\wh{\Gam})$. It is not even known if $\wh{\Gam}$ is torsion free when $\Gam$ has no torsion, hence $\wh{\Gam}$ could even have infinite cohomological dimension.
\end{rem}
%%%%%%%%%%%%%%%%%%%%

%%%%%%%%%%%%%%%%%%%%
\subsubsection*{Comments on the proofs and organization of the paper}
The proofs of all the main results of this paper rely on showing that certain classes in $H^*(\Gam_0 \bs \bbB^n, \bbZ)$ are pullbacks under a holomorphic map to some abelian variety $A$, namely the Albanese variety of the ball quotient. Goodness of $\pi_1(A)$ is then key in each proof. Postponing some technical results on arithmetic lattices of simplest type until \Cref{sec:Theta}, the proofs of \Cref{thm:Main}, \Cref{thm:MainTech}, and \Cref{thm:Main2} are in \Cref{sec:Proof}.
%%%%%%%%%%%%%%%%%%%%

%%%%%%%%%%%%%%%%%%%%
\subsubsection*{Acknowledgments}
This material is based upon work supported by Grants DMS-2203555 and DMS-2506896 from the National Science Foundation and support from the Simons Foundation [SFI-MPS-TSM-00014184, MS]. The author would additionally like to thank the Isaac Newton Institute for Mathematical Sciences, Cambridge, for support and hospitality during the programme `Discrete and profinite groups', where work on this paper was undertaken. This work was supported by EPSRC grant EP/Z000580/1. Finally, this paper was completed while the author was a CRM-Simons Professor at the Centre de Recherches Math\'ematiques, and he thanks both CRM and the Simons Foundation for support during that period.
%%%%%%%%%%%%%%%%%%%%

%%%%%%%%%%%%%%%%%%%%
\section{Proofs of the main results}\label{sec:Proof}
%%%%%%%%%%%%%%%%%%%%

See \cite{Reid1} for standard facts about profinite completions of discrete groups and \cite{Goldman} for basics on the geometry of $\bbB^n$ and its isometries. Since smooth compact ball quotients are aspherical, whenever $\Gam$ is cocompact and torsion-free this section freely employs the standard isomorphism between cohomology rings of $\Gam \bs \bbB^n$ and $\Gam$ with any given trivial coefficients.

The proofs of the primary theorems in this paper do not rely on the definition of an arithmetic lattice of simplest type, only some very special properties of their cohomology rings. To avoid complicating this section unnecessarily, that definition is postponed to \Cref{sec:Theta}. This section relies heavily on the following very deep but relatively simple-to-state result.

%%%%%%%%%%%%%%%%%%%%
\begin{thm}\label{thm:Wedge}
Suppose $\Gam < \PU(n,1)$ is a cocompact arithmetic lattice of simplest type.
\begin{enum}

\item There is a torsion-free subgroup $\Gam_0 \le \Gam$ of finite index so that, for all $\Del \le \Gam_0$ of finite index and $X = \Del \bs \bbB^n$, the image in $H^2(X, \bbC)$ of the first Chern class $K_X \in H^2(X, \bbZ)$ of the canonical line bundle on $X$ is in the image of the evaluation map
\begin{equation}\label{eq:Wedge2}
c_2 : \bigwedge\nolimits^2 H^1(X, \bbC) \lra H^2(X, \bbC)
\end{equation}
for the cup product.

\item Further, the subgroup $\Gam_0$ in Part 1.\ can be chosen so that the evaluation maps
\begin{equation}\label{eq:Wedge}
c_j : \bigwedge\nolimits^j H^1(X, \bbC) \lra H^j(X, \bbC)
\end{equation}
have nontrivial image for all $j \le 2n$ and all $X$ as in Part 1.

\item Suppose moreover that $\Gam$ is a congruence arithmetic lattice and that $j < \frac{n+1}{2}$. For any $r \in \bbN$, $\Gam_0$ can be chosen so any linearly independent $\phi_1, \dots, \phi_r$ in $H^j(X, \bbC)$ are contained in the image of $c_j$.

\end{enum}
\end{thm}
%%%%%%%%%%%%%%%%%%%%

Parts \emph{1}.\ and \emph{2}.\ follow fairly directly from results in the literature, so they are proved here. The proof of Part \emph{3}.\ is somewhat more involved, though special cases are known, so its proof is postponed until \Cref{sec:Theta}.

%%%%%%%%%%%%%%%%%%%%
\begin{pf}[Proof of 1.\ and 2.\ in \Cref{thm:Wedge}]
Let $\Gam$ be a lattice in $\PU(n,1)$. The theorem is well-known for $n = 1$, so assume $n \ge 2$. Note that residual finiteness of lattices in $\PU(n,1)$ (e.g., see \cite{StoverSurvey}) allows one to pass to a torsion-free subgroup of finite index. Moreover, standard differential topology implies that all properties asserted in the theorem pull back to any finite cover, so it suffices to find the initial lattice $\Gam_0$.

If $\Gam$ is a cocompact arithmetic lattice of simplest type, then Part \emph{1}.\ is a result proved jointly with Toledo \cite[Cor.\ 3.6]{StoverToledo2}. Let $X = \Gam_0 \bs \bbB^n$ satisfy Part \emph{1}. Given nontrivial $\al \in H^k(X, \bbC)$ and $\al^\prime \in H^{k^\prime}(X, \bbC)$ with $k + k^\prime \le n$, Venkataramana showed that there is a finite cover $p : Y \to X$ and an automorphism $g : Y \to Y$ (specifically, a \emph{Hecke correspondence}) so that
\begin{equation}\label{eq:Venky}
(g \circ p)^*(\al) \wedge p^*(\al^\prime) \neq 0
\end{equation}
in $H^{k + k^\prime}(Y, \bbC)$ \cite[Thm.\ 8]{Venky}. Strictly, Venkataramana proved this for $X$ a quotient by a congruence subgroup, but every arithmetic lattice is contained in a congruence lattice by \cite[Prop.\ 1.4]{BorelPrasad}, so one can pull back to an example commonly covered by the given lattice and a congruence lattice where \Cref{eq:Venky} holds to obtain the general conclusion. Induction on the case $k = k^\prime = 1$ implies that one can replace the initial $\Gam_0$ with a finite index subgroup such that $c_j$ has nontrivial image for all $j \le n$.

Set $X = \Gam_0 \bs \bbB^n$ and fix a nontrivial $\phi \in \im(c_j) \subseteq H^j(X, \bbC)$. To simplify notation, identify $K_X$ with its image in $H^2(X, \bbC)$ under tensor product. Smooth compact ball quotients are algebraic surfaces of general type with ample canonical bundle (e.g., see \cite[\S\S V.20 and VII.9]{BPV}). The hard Lefschetz theorem \cite[p.\ 122]{GriffithsHarris} then implies that wedge product with $K_X^{n - j}$ defines an isomorphism from $H^j(X, \bbC)$ to $H^{2n - j}(X, \bbC)$. Since $K_X$ is in the image of $c_2$, it follows that $\phi \wedge K_X^{n-j} \neq 0$ is in the image of $c_{2n - j}$. Thus $c_j$ has nontrivial image for all $j \le 2n$, which proves Part \emph{2}.
\end{pf}
%%%%%%%%%%%%%%%%%%%%

Let $X = \Gam_0 \bs \bbB^n$ be a compact quotient of $\bbB^n$ by a torsion-free lattice. Then there is an abelian variety $\mathrm{Alb}(X)$ called the \emph{Albanese variety} of $X$ and a holomorphic mapping ${\al : X \to \mathrm{Alb}(X)}$ so that
\[
\al^* : H^1(\mathrm{Alb}(X), \bbC) \lra H^1(X, \bbC)
\]
is an isomorphism. See \cite[p.\ 552]{GriffithsHarris}. Since wedge product commutes with pullback and the cohomology ring $H^*(\mathrm{Alb}(X), \bbC)$ is the exterior algebra on $H^1(\mathrm{Alb}(X), \bbC)$, \Cref{thm:Wedge} has the following immediate consequence.

%%%%%%%%%%%%%%%%%%%%
\begin{cor}\label{cor:Wedge}
Suppose $\Gam < \PU(n,1)$ is a cocompact arithmetic lattice of simplest type.
\begin{enum}

\item There is a torsion-free subgroup $\Gam_0 \le \Gam$ of finite index so that, for all $\Del \le \Gam_0$ of finite index and $X = \Del \bs \bbB^n$, the image in $H^2(X, \bbC)$ of the first Chern class $K_X \in H^2(X, \bbZ)$ of the canonical line bundle on $X$ is the pullback of a class $\psi_X \in H^2(\mathrm{Alb}(X), \bbC)$ under the Albanese map ${\al : X \to \mathrm{Alb}(X)}$.

\item Further, the subgroup $\Gam_0$ in Part 1.\ can be chosen so that the pullback maps
\begin{equation}\label{eq:Pullback}
\al^* : H^j(\mathrm{Alb}(X), \bbC) \lra H^j(X, \bbC)
\end{equation}
have nontrivial image for all $j \le 2n$ and all $X$ as in Part 1.

\item Suppose moreover that $\Gam$ is a congruence arithmetic lattice and that $j < \frac{n+1}{2}$. For any $r \in \bbN$, $\Gam_0$ can be chosen so any linearly independent $\phi_1, \dots, \phi_r$ in $H^j(X, \bbC)$ are contained in the image of $\al^*$.

\end{enum}
\end{cor}
%%%%%%%%%%%%%%%%%%%%

This is all the preparation needed to prove Theorems \ref{thm:Main} and \ref{thm:Main2}.

%%%%%%%%%%%%%%%%%%%%
\begin{pf}[Proof of Theorems \ref{thm:Main} and \ref{thm:Main2}]
Since Betti numbers only grow under finite covers, if the conclusion of \Cref{thm:Main} holds for $\Gam_0$ and is proved using reductions modulo $p$ of infinite order integral homology classes, then it holds for any finite index subgroup $\Del \le \Gam_0$. Thus it suffices to prove each theorem for $\Gam_0$. Fix a prime $p$ and set $X = \Gam_0 \bs \bbB^n$ for $\Gam_0$ the lattice provided by \Cref{cor:Wedge}. The Albanese map $\al : X \to \mathrm{Alb}(X)$ induces a commutative diagram
\[
\begin{tikzcd}
H^j(\wh{\bbZ}^d, \bbF_p) \arrow[r, "\wh{\al}_p^*"] \arrow[d, "\cong" left] & H^j(\wh{\Gam}_0, \bbF_p) \arrow[d, "\res_p"] \\
H^j(\bbZ^d, \bbF_p) \arrow[r, "\al_p^*"] & H^j(\Gam_0, \bbF_p) \\
H^j(\bbZ^d, \bbZ) \arrow[u, "r_p^\#"] \arrow[r, "\al^*"] & H^j(\Gam_0, \bbZ) \arrow[u, "r_p^\#" right]
\end{tikzcd}
\]
in cohomology, where $\bbZ^d = \pi_1(A)$ and $r_p^\#$ is the map induced by the reduction homomorphism $r_p : \bbZ \to \bbF_p$ on targets of cocycles. Note that $H^j(\wh{\bbZ}^d, \bbF_p) \cong H^j(\bbZ^d, \bbF_p)$ since $\bbZ^d$ is good. It follows that if $\phi \in H^j(\Gam_0, \bbZ)$ is $\al^*(\psi)$ for some $\psi$, then
\[
\wh{\phi}_p = \wh{\al}_p^*(\wh{\psi}_p)
\]
is an element of $H^j(\wh{\Gam}_0, \bbF_p)$ with $\res(\wh{\phi}_p) = r_p^\#(\phi)$, where $\wh{\psi}_p$ is identified with the image of $r_p^\#(\psi)$ in $H^j(\wh{\bbZ}^d, \bbF_p)$ under the restriction isomorphism.

It follows from \Cref{cor:Wedge} that the integral cup product map
\[
c^\bbZ_j : \bigwedge\nolimits^j H^1(\Gam_0, \bbZ) \to H^j(\Gam_0, \bbZ)
\]
contains an element of infinite order $\phi$ in its image. Indeed, if $c_j^\bbZ$ had image consisting only of torsion elements, then $c_j$ is the map induced by taking tensor products of $\bbZ$ cohomology groups with $\bbC$, and hence $c_j$ would have trivial image, contrary to hypothesis on $\Gam_0$. Then $r_p^\#(\phi) \in H^j(\Gam_0, \bbF_p)$ is nontrivial for all  but finitely many $p$. The previous paragraph then implies nontriviality of $H^j(\wh{\Gam}_0, \bbF_p)$ for all such $p$, as desired. This completes the proof of \Cref{thm:Main}.

\Cref{thm:Main2} follows in the same manner after replacing the single element $\phi$ with the span of $\phi_1, \dots, \phi_r$, which by assumption come from generators of a free submodule of $H^j(\Gam_0, \bbZ)$ of rank $r$. Fix $\psi_1, \dots, \psi_r$ in $H^j(\bbZ^n, \bbZ)$ with $\al^*(\psi_i) = \phi_i$ for all $i$ by Part \emph{3}.\ of \Cref{cor:Wedge}. The submodule $W_p$ of $H^j(\Gam_0, \bbF_p)$ spanned by $r_p^\#(\phi_1), \dots, r_p^\#(\phi_r)$ still has rank $r$ for all but finitely many $p$. Then one finds $\wh{\psi}_p$ in $H^j(\wh{\bbZ}^n, \bbF_p)$ as above so that $\wh{\phi}_i = \wh{\al}_p^*(\wh{\psi}_i)$ span the desired submodule of $H^j(\wh{\Gam}_0, \bbF_p)$ of rank $r$ over $\bbF_p$. This completes the proof.
\end{pf}
%%%%%%%%%%%%%%%%%%%%

Paying more attention to the kernel of the map from integral to complex cohomology groups leads to the proof of \Cref{thm:MainTech}.

%%%%%%%%%%%%%%%%%%%%
\begin{pf}[Proof of \Cref{thm:MainTech}]
Let $X$ be $\Gam_0 \bs \bbB^n$ for $\Gam_0 < \PU(n,1)$ the lattice provided by \Cref{cor:Wedge}. Then there is a class $\psi \in H^2(\mathrm{Alb}(X), \bbZ)$ so that $\al^*(\psi \otimes \bbC)$ equals the image of $K_X$ in $H^2(\Gam_0, \bbC)$ under the tensor product. Thus the pullback $\al_\bbZ^*$ on integral cohomology has
\[
\al_\bbZ^*(\psi) - K_X \in \ker\!\left(H^2(\Gam_0, \bbZ) \lra H^2(\Gam_0, \bbC) \right)
\]
and so $\al_\bbZ^*(\psi) = K_X + \sig$ for some torsion class $\sig \in H^2(\Gam_0, \bbZ)$. Then the $i^{th}$ exterior power $(K_X + \sig)^i$ is nontrivial in $H^{2i}(\Gam_0, \bbZ)$ for all $i \le n$, since the same is true for $K_X^i$. The main argument in the proof of \Cref{thm:Main} then implies that
\[
0 \neq (K_X + \sig)^i \!\!\!\!\!\pmod{p} \in \im\!\left(\mathrm{res} : H^{2i}(\wh{\Gam}_0, \bbF_p) \lra H^{2i}(\Gam_0, \bbF_p)\right)
\]
for all $i \le n$ and all but finitely many primes $p$. This proves the theorem for $\Gam_0$, and the result follows for any $\Del \le \Gam_0$ of finite index, albeit possibly with different exceptional primes, by pullback from $\Gam_0 \bs \bbB^n$. Thus the proof of the theorem is complete.
\end{pf}
%%%%%%%%%%%%%%%%%%%%

It now remains to justify Part \emph{3}.\ of \Cref{thm:Wedge}, which is the goal of the next section.

%%%%%%%%%%%%%%%%%%%%
\section{The theta correspondence}\label{sec:Theta}
%%%%%%%%%%%%%%%%%%%%

This section adapts arguments due to Bergeron, Millson, and Moeglin \cite{BMM} to prove Part \emph{3}.\ of \Cref{thm:Wedge}. Notation generally follows their paper; see \cite[\S 6]{BMM} and also \cite[\S 3]{StoverToledo2}.

To start, cocompact lattices in $\PU(n,1)$ of simplest type are constructed as follows. Let $F$ be a totally real number field of degree $d$ over $\bbQ$ and $E / F$ be a totally complex extension. If $(V, h)$ is a nondegenerate hermitian vector space over $E$ of dimension $n+1$, then the unitary group $\U(V)$ is a reductive algebraic group defined over $F$. Let $G_1$ denote the restriction of scalars of $\U(V)$ from $F$ to $\bbQ$. Assuming that the signature of $h$ is $(n,1)$ at exactly one real place of $F$ and that it is definite at all other real places, one obtains an isomorphism
\[
G_1(\bbR) \cong \U(n,1) \times \U(n+1)^{d-1}
\]
of real algebraic groups. A choice of maximal compact subgroup $K_\infty$ of $G_1(\bbR)$ identifies $\bbB^n$ with $G_1(\bbR) / K_\infty$.

Let $G$ be the restriction from $F$ to $\bbQ$ of the group of unitary similitudes of $(V, h)$. If $\bbA_\bbQ$ is the adeles of $\bbQ$, $\bbA_\bbQ^f$ is the finite adeles, and $K = \prod_p K_p$ an open compact subgroup of $G(\bbA_\bbQ^f)$ consisting of a product of open compact subgroups $K_p$ of $G(\bbQ_p)$, then one obtains a Shimura variety
\[
S(K) = G(\bbQ) \bs\!\left(\bbB^n \times G(\bbA_\bbQ^f)\right)\!/K
\]
that is a disjoint union of connected components $S(\Gam_j) = \Gam_j \bs \bbB^n$ \cite[\S 6.4]{BMM}. Any such $\Gam_j$ is a \emph{congruence arithmetic lattice of simplest type}, and any lattice ${\Gam < \PU(n,1)}$ commensurable with the image of $\Gam_j$ under the standard projection $G(\bbR) \to \PU(n,1)$ will be called an \emph{arithmetic lattice of simplest type}, and $\Gam$ is cocompact if and only if $F \neq \bbQ$.

The key for this section is determining when cohomology arises from \emph{theta lifting} \cite[\S 7]{BMM}, which is a powerful tool allowing one to construct cohomology classes and understand when they are in the image of the cup product map from smaller degrees. The following important result regarding the Hodge decomposition $H^{p,q}(X)$ of a smooth compact arithmetic ball quotient $X$ of simplest type will prove fundamental.

%%%%%%%%%%%%%%%%%%%%
\begin{thm}[Cor.\ 7.3 \cite{BMM}]\label{thm:BMMrange}
With the notation established in this section, suppose that $S$ is a connected component of $S(K)$ and $q \in \bbN$. If
\begin{equation}\label{eq:BMMrange}
3(a+b) + |a-b| < 2(n+1)
\end{equation}
then $H^{aq, bq}(S)$ is generated by classes of theta lifts.
\end{thm}
%%%%%%%%%%%%%%%%%%%%

Taking $q = 1$ and ensuring \Cref{eq:BMMrange} holds for all partitions $j = a+b$ with $a,b \ge 0$, one obtains the following corollary whose relevance to Part \emph{3}.\ of \Cref{thm:Wedge} is hopefully fairly evident.

%%%%%%%%%%%%%%%%%%%%
\begin{cor}\label{cor:BMMrangeAll}
With $S$ as in \Cref{thm:BMMrange}, $H^j(S)$ is generated by classes of theta lifts for all $j < \frac{n+1}{2}$.
\end{cor}
%%%%%%%%%%%%%%%%%%%%

The following is the essential result of this section.

%%%%%%%%%%%%%%%%%%%%
\begin{thm}\label{thm:MyCup}
Following the notation established in this section, suppose $S$ is a component of $S(K)$, $j < \frac{n+1}{2}$, and $\phi \in H^j(S, \bbC)$. Then there is a congruence cover $p : S^\prime \to S$ so that $p^*(\phi)$ is in the image of the evaluation map
\[
\bigwedge\nolimits^j H^1(S^\prime, \bbC) \lra H^j(S^\prime, \bbC)
\]
for the cup product.
\end{thm}
%%%%%%%%%%%%%%%%%%%%

%%%%%%%%%%%%%%%%%%%%
\begin{pf}
Set $Sh(G) = G(\bbQ) \bs\!\left(\bbB^n \times G(\bbA_\bbQ^f)\right)\!$, which has cohomology ring
\[
H^*(Sh(G), \bbC) = \varinjlim_L H^*(S(L), \bbC)
\]
from the limit over all open compact subgroups $L < G(\bbA_\bbQ^f)$. Given any ${\phi \in H^j(S, \bbC)}$, let $\wt{\phi} \in H^j(Sh(G), \bbC)$ be the associated limiting class, which is nontrivial since $\phi$ pulls back to a nontrivial class on any finite cover.

Note that $\phi$ is generated by theta lifts by \Cref{cor:BMMrangeAll}. As described in \cite[\S 9.2]{BMM}, this means more precisely that $H^{a,b}(Sh(G), \bbC)$ is generated by certain automorphic functions $\thet^f_{\psi, \chi, \varphi_{a,b} \otimes \phi_f}$ constructed from special cocycles
\[
\psi_{a,b} \in \Hom_{K_\infty}\!\left(\bigwedge\nolimits^{a,b} \frakp, \calP(V \otimes W)\right)
\]
for $\frakp$ the maximal parabolic subalgebra of the Lie algebra of $G(\bbR)$ and $\calP(V \otimes W)$ the Fock space for the oscillator representation defined in \cite[\S 4.2]{BMM}, where $W$ is the space defined in the proof of \cite[Thm.\ 9.3]{BMM}. Then Prop.\ 5.4, Prop.\ 5.19, and the calculation in Eq.\ (8.10) of \cite{BMM} combine exactly as in the proof of \cite[Thm.\ 9.3]{BMM} to imply that the special cocycles $\psi_{a,b}$ are iterated wedge products of the classes $\psi_{1,0}$ and $\psi_{0,1}$ and that this property remains true in $H^*(Sh(G), \bbC)$.

In particular, $\wt{\phi}$ is in the image of
\[
\wt{c}_j : \bigwedge\nolimits^j H^1(Sh(G), \bbC) \lra H^j(Sh(G), \bbC)
\]
and thus there is a congruence cover $p : S^\prime \to S$ so that
\[
p^*(\phi) \in \im\!\left(\bigwedge\nolimits^j H^1(S^\prime, \bbC) \lra H^j(S^\prime, \bbC) \right)
\]
exactly as in \cite[Thm.\ 3.4]{StoverToledo2}. Indeed, that the desired statement holds in the direct limit means that there is a compact open subgroup $K^\prime \le K$ so that the pullback $\phi^\prime$ of $\phi$ to some component $S^\prime$ of $S(K^\prime)$ has the property that $\phi^\prime$ is in the image of the cup product from $\bigwedge\nolimits^j H^1(S^\prime, \bbC)$. This completes the proof of the theorem.
\end{pf}
%%%%%%%%%%%%%%%%%%%%

This is the required preparation to complete the goal of this section.

%%%%%%%%%%%%%%%%%%%%
\begin{pf}[Proof of Part 3.\ of \Cref{thm:Wedge}]
Suppose $\Gam < \PU(n,1)$ is a cocompact congruence arithmetic lattice of simplest type. With notation as in this section, given $\phi \in H^j(\Gam, \bbC)$, pass to a torsion free congruence subgroup of finite index associated with an open compact subgroup $K$ of the associated algebraic group $G(\bbA_\bbQ^f)$ so that $\phi$ pulls back to a class (still denoted $\phi$) in $H^j(S, \bbC)$ for some component $S$ of $S(K)$. \Cref{thm:MyCup} then supplies a congruence cover $p : S^\prime \to S$ so that $p^*(\phi)$ is in the image of the cup product from $H^1(S^\prime, \bbC)$. Repeating over a finite collection of linearly independent elements of $H^1(\Gam, \bbC)$ and pullback to a common congruence cover proves the desired result.
\end{pf}
%%%%%%%%%%%%%%%%%%%%

%%%%%%%%%%%%%%%%%%%%
\begin{rem}
Using the hard Lefschetz theorem as in the proof of Part \emph{2}.\ of \Cref{thm:Wedge}, one can also obtain some surjectivity of the restriction map in higher degrees.
\end{rem}
%%%%%%%%%%%%%%%%%%%%

%%%%%%%%%%%%%%%%%%%%
\bibliography{ProfiniteCanonical}

\begin{thebibliography}{10}

\bibitem{BFMS2}
Uri Bader, David Fisher, Nicholas Miller, and Matthew Stover.
\newblock Arithmeticity, superrigidity and totally geodesic submanifolds of
  complex hyperbolic manifolds.
\newblock {\em Invent. Math.}, 233(1):169--222, 2023.
\newblock URL: \url{https://doi.org/10.1007/s00222-023-01186-5}.

\bibitem{BU}
Gregorio Baldi and Emmanuel Ullmo.
\newblock Special subvarieties of non-arithmetic ball quotients and {H}odge
  theory.
\newblock {\em Ann. of Math. (2)}, 197(1):159--220, 2023.
\newblock URL: \url{https://doi.org/10.4007/annals.2023.197.1.3}.

\bibitem{BPV}
W.~Barth, C.~Peters, and A.~Van~de Ven.
\newblock {\em Compact complex surfaces}, volume~4 of {\em Ergebnisse der
  Mathematik und ihrer Grenzgebiete (3)}.
\newblock Springer-Verlag, 1984.
\newblock URL: \url{https://doi.org/10.1007/978-3-642-96754-2}.

\bibitem{BHH}
Gottfried Barthel, Friedrich Hirzebruch, and Thomas H\"{o}fer.
\newblock {\em Geradenkonfigurationen und {A}lgebraische {F}l\"{a}chen}.
\newblock Aspects of Mathematics, D4. Friedr. Vieweg \& Sohn, Braunschweig,
  1987.
\newblock URL: \url{https://doi.org/10.1007/978-3-322-92886-3}.

\bibitem{BMM}
Nicolas Bergeron, John Millson, and Colette Moeglin.
\newblock The {H}odge conjecture and arithmetic quotients of complex balls.
\newblock {\em Acta Math.}, 216(1):1--125, 2016.
\newblock URL: \url{https://doi.org/10.1007/s11511-016-0136-2}.

\bibitem{BorelPrasad}
Armand Borel and Gopal Prasad.
\newblock Finiteness theorems for discrete subgroups of bounded covolume in
  semi-simple groups.
\newblock {\em Inst. Hautes \'{E}tudes Sci. Publ. Math.}, (69):119--171, 1989.
\newblock URL: \url{http://www.numdam.org/item?id=PMIHES_1989__69__119_0}.

\bibitem{Clozel}
Laurent Clozel.
\newblock On the cohomology of {K}ottwitz's arithmetic varieties.
\newblock {\em Duke Math. J.}, 72(3):757--795, 1993.
\newblock URL: \url{https://doi.org/10.1215/S0012-7094-93-07229-8}.

\bibitem{DerauxForget}
Martin Deraux.
\newblock Forgetful maps between {D}eligne-{M}ostow ball quotients.
\newblock {\em Geom. Dedicata}, 150:377--389, 2011.
\newblock URL: \url{https://doi.org/10.1007/s10711-010-9511-x}.

\bibitem{DiCerboStover2}
Luca~F. Di~Cerbo and Matthew Stover.
\newblock Bielliptic ball quotient compactifications and lattices in {$\rm
  PU(2,1)$} with finitely generated commutator subgroup.
\newblock {\em Ann. Inst. Fourier (Grenoble)}, 67(1):315--328, 2017.
\newblock URL: \url{https://doi.org/10.5802/aif.3083}.

\bibitem{DiCerboStover1}
Luca~F. Di~Cerbo and Matthew Stover.
\newblock Classification and arithmeticity of toroidal compactifications with
  {$3\overline c_2=\overline c_1^2=3$}.
\newblock {\em Geom. Topol.}, 22(4):2465--2510, 2018.
\newblock URL: \url{https://doi.org/10.2140/gt.2018.22.2465}.

\bibitem{DiCerboStover3}
Luca~F. Di~Cerbo and Matthew Stover.
\newblock Punctured spheres in complex hyperbolic surfaces and bielliptic ball
  quotient compactifications.
\newblock {\em Trans. Amer. Math. Soc.}, 372(7):4627--4646, 2019.
\newblock URL: \url{https://doi.org/10.1090/tran/7650}.

\bibitem{Goldman}
William~M. Goldman.
\newblock {\em Complex hyperbolic geometry}.
\newblock Oxford Mathematical Monographs. The Clarendon Press, 1999.
\newblock Oxford Science Publications.

\bibitem{GriffithsHarris}
Phillip Griffiths and Joseph Harris.
\newblock {\em Principles of algebraic geometry}.
\newblock Wiley Classics Library. John Wiley \& Sons, Inc., 1994.
\newblock URL: \url{https://doi.org/10.1002/9781118032527}.

\bibitem{GJZZ}
F.~Grunewald, A.~Jaikin-Zapirain, and P.~A. Zalesskii.
\newblock Cohomological goodness and the profinite completion of {B}ianchi
  groups.
\newblock {\em Duke Math. J.}, 144(1):53--72, 2008.
\newblock URL: \url{https://doi.org/10.1215/00127094-2008-031}.

\bibitem{Hirzebruch}
F.~Hirzebruch.
\newblock Chern numbers of algebraic surfaces: an example.
\newblock {\em Math. Ann.}, 266(3):351--356, 1984.
\newblock URL: \url{https://doi.org/10.1007/BF01475584}.

\bibitem{Kapovich}
Michael Kapovich.
\newblock On normal subgroups of the fundamental groups of complex surfaces.
\newblock Preprint, 1998.
\newblock URL: \url{https://www.math.ucdavis.edu/~kapovich/EPR/koda.pdf}.

\bibitem{Kazhdan}
David Kazhdan.
\newblock Some applications of the {W}eil representation.
\newblock {\em J. Analyse Math.}, 32:235--248, 1977.
\newblock URL: \url{https://doi.org/10.1007/bf02803582}.

\bibitem{KRS}
Alexander Kolpakov, Alan~W. Reid, and Leone Slavich.
\newblock Embedding arithmetic hyperbolic manifolds.
\newblock {\em Math. Res. Lett.}, 25(4):1305--1328, 2018.
\newblock URL: \url{https://doi.org/10.4310/MRL.2018.v25.n4.a12}.

\bibitem{Livne}
Ron~Aharon Livne.
\newblock {\em On certain covers of the universal elliptic curve}.
\newblock ProQuest LLC, Ann Arbor, MI, 1981.
\newblock Thesis (Ph.D.)--Harvard University.
\newblock URL:
  \url{http://gateway.proquest.com/openurl?url_ver=Z39.88-2004&rft_val_fmt=info:ofi/fmt:kev:mtx:dissertation&res_dat=xri:pqdiss&rft_dat=xri:pqdiss:0372704}.

\bibitem{Lochak}
Pierre Lochak.
\newblock Results and conjectures in profinite {T}eichm\"{u}ller theory.
\newblock In {\em Galois-{T}eichm\"{u}ller theory and arithmetic geometry},
  volume~63 of {\em Adv. Stud. Pure Math.}, pages 263--335. Math. Soc. Japan,
  Tokyo, 2012.
\newblock URL: \url{https://doi.org/10.2969/aspm/06310263}.

\bibitem{Reid1}
Alan~W. Reid.
\newblock Profinite properties of discrete groups.
\newblock In {\em Groups {S}t {A}ndrews 2013}, volume 422 of {\em London Math.
  Soc. Lecture Note Ser.}, pages 73--104. Cambridge Univ. Press, 2015.

\bibitem{Serre}
Jean-Pierre Serre.
\newblock {\em Cohomologie galoisienne}, volume~5 of {\em Lecture Notes in
  Mathematics}.
\newblock Springer-Verlag, Berlin, fourth edition, 1994.
\newblock URL: \url{https://doi.org/10.1007/BFb0108758}.

\bibitem{Stix}
Jakob Stix.
\newblock {\em Projective anabelian curves in positive characteristic and
  descent theory for log-\'{e}tale covers}, volume 354 of {\em Bonner
  Mathematische Schriften}.
\newblock Universit\"{a}t Bonn, Mathematisches Institut, 2002.
\newblock Dissertation, Rheinische Friedrich-Wilhelms-Universit\"{a}t Bonn.

\bibitem{StoverFake}
Matthew Stover.
\newblock Algebraic fundamental groups of fake projective planes.
\newblock To appear in J. Eur. Math. Soc. (JEMS).
\newblock URL: \url{https://arxiv.org/abs/2205.13991}.

\bibitem{StoverExs}
Matthew Stover.
\newblock Products of curves as ball quotients.
\newblock Submitted.
\newblock URL: \url{https://arxiv.org/abs/2312.05699}.

\bibitem{StoverSurvey}
Matthew Stover.
\newblock Residual finiteness and discrete subgroups of {L}ie groups.
\newblock To appear in Zariski Dense Subgroups, Number Theory and Geometric
  Applications.
\newblock URL: \url{https://arxiv.org/abs/2407.07680}.

\bibitem{StoverProfinite}
Matthew Stover.
\newblock Lattices in {${\rm PU}(n,1)$} that are not profinitely rigid.
\newblock {\em Proc. Amer. Math. Soc.}, 147(12):5055--5062, 2019.
\newblock URL: \url{https://doi.org/10.1090/proc/14763}.

\bibitem{StoverToledo2}
Matthew Stover and Domingo Toledo.
\newblock Residual finiteness for central extensions of lattices in {${\rm
  PU}(n, 1)$} and negatively curved projective varieties.
\newblock {\em Pure Appl. Math. Q.}, 18(4):1771--1797, 2022.
\newblock URL: \url{https://doi.org/10.4310/pamq.2022.v18.n4.a15}.

\bibitem{StoverToledo1}
Matthew Stover and Domingo Toledo.
\newblock Residually finite lattices in {$\widetilde{PU(2,1)}$} and fundamental
  groups of smooth projective surfaces.
\newblock {\em Michigan Math. J.}, 72:559--597, 2022.
\newblock URL: \url{https://doi.org/10.1307/mmj/20217215}.

\bibitem{ToledoMap}
Domingo Toledo.
\newblock Maps between complex hyperbolic surfaces.
\newblock volume~97, pages 115--128. 2003.
\newblock Special volume dedicated to the memory of Hanna Miriam Sandler
  (1960--1999).
\newblock URL: \url{https://doi.org/10.1023/A:1023691505890}.

\bibitem{Venky}
T.~N. Venkataramana.
\newblock Cohomology of compact locally symmetric spaces.
\newblock {\em Compositio Math.}, 125(2):221--253, 2001.
\newblock URL: \url{https://doi.org/10.1023/A:1002600432171}.

\end{thebibliography}
%%%%%%%%%%%%%%%%%%%%

%%%%%%%%%%%%%%%%%%%%
\end{document}